\documentclass{conm-p-l}
\usepackage{amsmath,amsfonts,amscd,amsthm,latexsym,amsxtra,amssymb}

{
}

\newcommand{\enquote}[1]{``#1''}

\begin{document}

\newcommand{\comment}[1]{{\bf #1}}
\newcommand{\forget}[1]{}

\newcommand{\mathcala}{{\mathcal{A}}}
\newcommand{\mathcalc}{{\mathcal C}}
\newcommand{\mathcalf}{{\mathcal F}}
\newcommand{\mathcalm}{{\mathcal M}}
\newcommand{\mathcaln}{{\mathcal N}}

\newcommand{\rr}{{\mathbb R}}
\newcommand{\zz}{{\mathbb Z}}
\newcommand{\nn}{{\mathbb N}}
\newcommand{\zmod}{\zz/2}
\newcommand{\cc}{{\mathbb C}}
\newcommand{\qq}{{\mathbb Q}}
\newcommand{\onto}{\twoheadrightarrow}
\newcommand{\map}{\operatorname{map}}
\newcommand{\colim}{\operatorname{colim}}
\newcommand{\tr}{\operatorname{tr}}
\newcommand{\pr}{\operatorname{pr}}
\newcommand{\id}{\operatorname{id}}
\newcommand{\im}{\operatorname{im}}
\newcommand{\coker}{\operatorname{coker}}
\newcommand{\clos}{\operatorname{clos}}
\newcommand{\subgroup}{<}
\newcommand{\tensor}{\otimes}
\newcommand{\defi}{\operatorname{def}}
\newcommand{\cent}{\operatorname{cent}}
\newcommand{\cone}{\operatorname{cone}}
\newcommand{\cyl}{\operatorname{cyl}}
\newcommand{\aut}{\operatorname{aut}}
\newcommand{\res}{\operatorname{res}}
\newcommand{\sing}{\operatorname{sing}}
\newcommand{\ind}{\operatorname{ind}}
\newcommand{\capac}{c}
\newcommand{\measurable}{measurable}
\newcommand{\cofinalm}{cofinal-measurable}

\newtheorem{theorem}{Theorem}[section]
\newtheorem{lemma}[theorem]{Lemma}
\newtheorem{corollary}[theorem]{Corollary}
\newtheorem{proposition}[theorem]{Proposition}

\theoremstyle{definition}
\newtheorem{definition}[theorem]{Definition}
\newtheorem{example}[theorem]{Example}
\newtheorem{notation}[theorem]{Notation}
\newtheorem{convention}[theorem]{Convention}
\newtheorem{remark}[theorem]{Remark}

{\catcode`@=11\global\let\c@equation=\c@theorem}
\renewcommand{\theequation}{\thetheorem}


\newcommand{\comsquare}[8]{
\begin{center}
$\begin{CD}
#1 @>#2>> #3\\
@V{#4}VV @VV{#5}V\\
#6 @>>#7> #8
\end{CD}$
\end{center}}

\newcommand{\textshortexactsequence}[5]{\mbox{$0\longrightarrow
#1\stackrel{#2}{\longrightarrow} #3\stackrel{#4}{\longrightarrow}
#5\longrightarrow 0$}}

\newcommand{\squarematrix}[4]
{\left( \begin{array}{cc} #1 & #2 \\ #3 & #4\end{array} \right)}



\typeout{--------------------- Introduction  -----------------------}

\title[Novikov-Shubin invariants for arbitrary group
actions]{Novikov-Shubin invariants for arbitrary group actions and
  their positivity}

\author{Wolfgang L\"uck}
\author{Holger Reich}
\author{Thomas Schick}
\address{Fachbereich Mathematik ---
Universit\"at M\"unster\\
Einsteinstr.~62 ---
         48149 M\"unster }

\email{wolfgang.lueck@math.uni-muenster.de}
\email{holger.reich@math.uni-muenster.de}
\email{thomas.schick@math.uni-muenster.de}

\dedicatory{Dedicated to Mel Rothenberg on the occasion of his 65th birthday}
\urladdr{http://www.uni-muenster.de/math/u/lueck/}

\begin{abstract}
We extend the notion of Novikov-Shubin invariant for
free $\Gamma$-$CW$-complexes of finite type to
spaces with arbitrary $\Gamma$-actions and prove some statements
about their positivity. In particular we apply this to classifying
spaces of discrete groups.
\end{abstract}

\maketitle
\section{Introduction}

In \cite{Lueck (1998a)} the first author
extended the notion of $L^2$-Betti number for
free $\Gamma$-$CW$-complexes of finite
type to topological spaces
with arbitrary $\Gamma$-actions. The key
 ingredient there is the notion of
a dimension function for arbitrary modules
over a finite von Neumann algebra $\mathcala$,
extending the classical notion for finitely
generated projective $\mathcala$-modules
defined in terms of the von Neumann trace of an associated projection.
These notions turned out to be useful in
particular in the case where $\Gamma$ is
amenable (see \cite{Lueck (1998a)},
\cite{Lueck (1998b)}).

In this paper we carry out an analogous
program for Novikov-Shubin invariants.
So we will introduce for arbitrary $\mathcala$-modules in Section \ref{Capacity of modules}
the equivalent notion of capacity (which is essentially
the inverse of the Novikov-Shubin invariant
and was introduced for finitely presented
$\mathcala$-modules in \cite{Farber (1996)} and \cite{Lueck (1997a)})
 and study its main properties.
This enables us to define
the $p$-th capacity for a space $X$ with an
arbitrary action of a discrete group. Originally Novikov-Shubin
invariants were defined in terms of the heat kernel of the universal
covering of a compact Riemannian manifold in
\cite{Novikov-Shubin (1987)}, \cite{Novikov-Shubin (1987b)}.

We will use the extension to study Novikov-Shubin invariants of groups in Section
\ref{Capacity of groups}. The key observation is that a group $\Gamma$
which may have a model of finite type for
$B\Gamma$ may contain interesting normal
subgroups for which the classical definition does not apply
because its classifying space is not even of finite type.
We are in particular interested
in the conjecture that for a regular covering
of a $CW$-complex of finite type
the Novikov-Shubin invariants are always positive,
or equivalently, the capacities are
always finite
\cite[Conjecture 7.1]{Lott-Lueck (1995)}.
Our main results in this direction
are Theorems
\ref{estimates on capac, extensions and local information} and
 \ref{estimate on capac_1 of an elementary amenable group}
which say among other things
that $\alpha_p(\Gamma) \ge 1$ for all
$p \ge 1$ if $\Gamma$ contains $\zz^n$
as normal subgroup for some $n \ge 1$
and that $\alpha_p(\Gamma) \ge 1$ for $p = 1,2$
if $\Gamma$ is elementary amenable and
contains no infinite locally finite subgroup.
In particular the first and second
Novikov-Shubin invariants of the universal covering of a $CW$-complex $X$
of finite type are greater or equal to $1$
if $\pi_1(X)$ satisfies the condition above
since they agree with the ones for $\pi_1(X)$.

We mention that $\alpha_p$ of a space or group
is the Novikov-Shubin invariant associated
to the $p$-th differential. Sometimes Novikov-Shubin invariants
are also defined in terms of the Laplacian. If we denote the latter
ones by $\widetilde{\alpha}_p$, the connection between these two invariants
is $\widetilde{\alpha}_p = \min\{\alpha_{p+1},\alpha_p\}$. Our normalization
is such that $\alpha_1 = \widetilde{\alpha}_0$ of the universal covering
of $S^1$ is $1$. The $p$-th capacity will essentially be the inverse
of the $(p+1)$-th Novikov-Shubin invariant.



\typeout{--------------------- section 1 -----------------------}
\section{Capacity of modules}
\label{Capacity of modules}

There is the notion of the Novikov-Shubin invariant $\alpha(M) \in \left[0,\infty \right] \amalg \{ \infty^{+} \} $
of a finitely presented $\mathcala$-module $M$ as defined in
\cite[Definition 3.1]{Lueck (1997a)}.
It is the Novikov-Shubin invariant of the spectral density function
of a morphism of Hilbert $\mathcala$-modules
$f: l^2(\mathcala)^m \to l^2(\mathcala)^n$
which corresponds to a presentation
matrix $A \in M(m,n,\mathcala)$ for $M=\coker(A: \mathcala^m \to \mathcala^n)$.
We want to extend it to arbitrary $\mathcala$-modules. For convenience
we will use the  notion of capacity (see \cite[4.8]{Farber (1996)})
which is essentially the
reciprocal of the Novikov-Shubin invariant and extend it to arbitrary $\mathcala$-modules.
In the sequel we will use the notion and properties
of the dimension $\dim(M)$ of an $\mathcala$-module
introduced in \cite[Theorem 0.6]{Lueck (1997a)}.

\begin{definition}
\label{definition of special and cofinal-special module}
We call an $\mathcala$-module $M$ {\em {\measurable}}
if it is the quotient of a finitely presented $\mathcala$-module
$L$ with $\dim(L) = 0$.
We call an $\mathcala$-module $M$ {\em {\cofinalm}}
if each finitely generated
submodule is {\measurable}. In particular this implies $\dim M=0$.\qed
\end{definition}

We will see later
that the  following definition of capacity is particularly well
behaved on the the class of {\cofinalm} modules.

\begin{definition}
\label{definition of capacity}
Let $L$ be a finitely presented $\mathcala$-module with
$\dim(L) = 0$.
Define its {\em capacity}
$$\capac(L) \in \{0^-\}\cup [0,\infty] $$
by
$$\capac(L) := \frac{1}{\alpha(L)},$$
where $0^-$ is a new formal symbol different from $0$, and
$0^{-1}=\infty$, $\infty^{-1}=0$ and $(\infty^+)^{-1}=0^-$.
Let $M$ be a {\measurable} $\mathcala$-module. Define
$$\capac'(M) ~ := ~
\inf\{\capac(L) \mid L \mbox{ finitely presented},~
\dim(L) = 0,~  M \mbox{ quotient of } L\}.$$
Let $N$ be an arbitrary $\mathcala$-module. Define
$$\capac''(N) ~ := ~ \sup\{\capac'(M) \mid M \mbox{ {\measurable}},~
M \subset N\}. \qed $$
\end{definition}
Note that $\dim(N)$ is not necessarily  zero. In fact
$c(N)$ measures the size of the largest zero-dimensional
submodule of $N$ (compare
\cite[2.15]{Lueck (1998a)}).

The invariants take value in
$\{0^-\} \amalg [0,\infty]$. We define an order $ < $ on this set
by the usual one on
$[0,\infty)$ and the rule
$0^- < r < \infty$ for $r \in [0,\infty)$.
For two elements $r,s\in [0^-,\infty]$ we define another element
$r+s=s+r$ in this set by
the ordinary addition in $[0,\infty)$ and by
$$\begin{array}{ccccccl}
0^- + r & = & r & \text{and} & \infty + r & = &
\infty \quad\text{for } r \in [0^-,\infty]
\end{array}$$

We have introduced $\capac(M)$ as the inverse of $\alpha(M)$
because then $\capac(M)$ becomes bigger if $M$ becomes bigger
and some of the formulas become nicer.
Notice that for a finitely presented
$\mathcala$-module $M$ with $\dim(M) =0$ we have
$\capac(M) = 0^-$ if and only if $M$ is trivial.
Hence a {\measurable} $\mathcala$-module $M$ satisfies $\capac'(M) = 0^-$
if and only if it is trivial.
A {\measurable} $\mathcala$-module is finitely
generated but not necessarily finitely presented.
We have for an arbitrary $\mathcala$-module $N$ that $\capac''(N) = 0^-$
if and only if any $\mathcala$-map
$f: M \longrightarrow N$
from a finitely presented $\mathcala$-module
$M$ with $\dim(M) = 0$ to $N$ is trivial.
A {\cofinalm} $\mathcala$-module $M$ is trivial if and only if
$\capac''(M) =0^- $ (see also Example
\ref{finitely generated module with dim = 0 and kappa = 0^-}).

We will show that our different definitions of capacity coincide and
prove their basic properties. For this we need
\begin{lemma}
  \label{exact_seq_fp}
  Let $0\to K\to P\to Q\to 0$ be an exact sequence of finitely
  presented $\mathcala$-modules of dimension zero. Then
  \begin{equation*}
    \begin{split}
      \capac(K)&\le \capac(P),\qquad\capac(Q)\le \capac(P)\\
      \capac(P)&\le \capac(K)+\capac(Q).
  \end{split}
\end{equation*}
\end{lemma}
\begin{proof}
  By \cite[3.4]{Lueck (1997a)} we find resolutions $0\to F_K\to F_K\to
  K\to 0$ and $0\to F_Q\to F_Q\to Q\to 0$ with $F_K$ and $F_Q$
  finitely generated free. Now we can construct a resolution $0\to
  F\to F\to P\to 0$ with $F=F_K\oplus F_Q$ which fits into  a short
  exact sequence of resolutions. Application of \cite[Lemma
  1.12]{Lott-Lueck (1995)} to this situation together
  with the equivalence of  finitely generated Hilbert
  $\mathcala$-modules and finitely generated projective algebraic $\mathcala$-modules as in
  \cite{Lueck (1997a)} gives the result.
\end{proof}

\begin{proposition}
\label{capac well def}
\begin{enumerate}
\item \label{capac = capac'}
If $M$ is a finitely presented
$\mathcala$-module which satisfies $\dim (M)=0$,
then $M$ is {\measurable} and
$$\capac'(M) ~ = ~ \capac (M).$$
\item \label{capac' = capac''}
If $M$ is a {\measurable} $\mathcala$-module,
then $M$ is {\cofinalm} and
$$\capac''(M) ~ = ~ \capac'(M).$$
\end{enumerate}
\end{proposition}

\begin{notation}
\label{kappa = kappa' = kappa'' dim^{prime}}
In view of Proposition
\ref{capac well def}
we will not distinguish between $\capac$, $\capac'$ and $\capac''$
in the sequel.
\end{notation}
\begin{proof}[Proof of Proposition \ref{capac well def}]
  If $M$ and $N$ are finitely presented with dimension zero and
  $p:N\to M$ is surjective, then $c(M)\le c(N)$ by
  Lemma \ref{exact_seq_fp}. We use that $\ker p$ is finitely
  presented by \cite[Theorem 0.2]{Lueck (1997a)}. The
  first statement follows.

  Suppose $M$ is {\measurable} and $M_0\subset M$ is finitely
  generated. We have to show that $M_0$ is {\measurable} and $\capac'(M_0)\le
  \capac'(M)$. Suppose $L$ is finitely presented with
$\dim(L)=0$ and there is an epimorphism
$f: L \longrightarrow  M$.
As $M_0$ is finitely generated, we can find a finitely generated module
$K \subset L$
with $f(K) = M_0$. As $K$ is finitely generated and
$L$ is finitely presented, $L/K$ is finitely presented.
Since the category of finitely presented $\mathcala$-modules is
abelian \cite[Theorem 0.2]{Lueck (1997a)},
$K$ is finitely presented. Hence $M_0$ is {\measurable}. Moreover
$c(K)\le c(L)$ by   Lemma \ref{exact_seq_fp}. Therefore
\begin{equation*}
  \capac'(M_0)\le \inf_{K \text{ as above}}\capac(K)\le \inf_{L\text{
      as above}}\capac(L)=\capac'(M).
\end{equation*}
\end{proof}

Before we give a list of the basic properties of capacity, we state a
simple lemma which will be used repeatedly during the proof.
\begin{lemma}\label{projinf}
  Suppose $M$ is a {\measurable} $\mathcala$-module and $p:F\to M$ a
  projection of a finitely generated free module onto $M$. Then
  \begin{equation*}
     \capac(M)= \inf\{\capac(F/K)|\; K\subset \ker p\text{ finitely
       generated with }\dim(F/K)=0 \} .
   \end{equation*}
   The set on the right hand side is  nonempty if and only if $M$ is \measurable.
\end{lemma}
\begin{proof}
 Let $d$ denote the number on the right. Every $F/K$ as above is
finitely presented and projects onto $M$. From Definition
\ref{definition of capacity} we get $d \geq c(M)$. On the other hand let
$q:L \to M$ be an epimorphism with $L$ finitely presented and $\dim L = 0$.
We can lift $p$ to a map $f$ with $q \circ f = p$.
Since the category of finitely presented modules is abelian \cite[Theorem 0.2]{Lueck (1997a)}
the kernel of $f$ is finitely generated. Moreover $\dim(F/ \ker f ) = \dim \im f \leq \dim L =0$,
$\ker f \subset \ker p$ and Lemma \ref{exact_seq_fp} implies $c(F/ \ker f ) \leq c(L)$.
So for every $L$ we found an $F/K$ with $c(F/K) \leq c(L)$. This implies
$d \leq c(M)$.
\end{proof}

\begin{theorem}
\label{basic properties of capacity}
\begin{enumerate}
\item \label{capacity and exact sequences}
Let
\mbox{$0 \longrightarrow M_0 \stackrel{i}{\longrightarrow} M_1
  \stackrel{p}{\longrightarrow}
M_2 \longrightarrow 0$} be an exact sequence of $\mathcala$-modules.
Then
\begin{enumerate}
\item\label{submod} $\capac(M_0) \le \capac(M_1)$.
\item\label{quotmod} $\capac(M_2) \le \capac(M_1)$, provided that $M_1$
is {\cofinalm}.
\item\label{extension} $\capac(M_1) \le \capac(M_0) + \capac(M_2)$ if
  $\dim M_1=0$.
\end{enumerate}

\item \label{unions}
Let $M=\bigcup_{i\in I} M_i$ be a directed union of submodules. Then
\begin{equation*}
  \capac(M) =    \sup\{\capac(M_i) \mid i \in I\}.
\end{equation*}

\item \label{capac and colimits}
Let $M$ be the colimit $\colim_{i \in I} M_i$ of a directed system
of $\mathcala$-modules with structure maps $\phi_{ij}:M_i\to M_j$.
Then
\begin{eqnarray*}
\capac(M) & \le & \liminf_{i\in I}\capac(M_i)\quad\left(:= \sup_{i\in
  I}\{\inf_{j\ge i}\capac(M_j)\} \right).
\end{eqnarray*}
If every $M_i$ is {\measurable} and $\phi_{ij}$ is
surjective for every $j\ge i$, then
\begin{equation*}
  \capac(M) = \inf_{i\in I}\capac(M_i).
\end{equation*}

\item \label{capac and direct sums}
Let $\{M_i \mid i \in I\}$ be a family of
$\mathcala$-modules.
Then
\begin{eqnarray*}
\capac\left(\oplus_{i \in I} M_i\right) & = &
\sup\{\capac(M_i) \mid i \in I\}.
\end{eqnarray*}
\end{enumerate}
\end{theorem}

\begin{remark}\label{basicdim}
  Because zero-dimensional modules will be most important for
  us, we give a reminder of basic properties of the dimension, which
  are stated in or follow  from \cite{Lueck (1998a)}:
  \begin{enumerate}
  \item If $0\to M_0\to M_1\to M_2\to 0$ is an exact sequence of
    $\mathcala$-modules, then
    \begin{equation*}
      \dim(M_1)=0 \quad\iff\quad \dim(M_0)=\dim(M_2)=0.
    \end{equation*}
  \item $\dim(\oplus_{i\in I}(M_i))=0$
    $\iff$ $\dim(M_i)=0$ for all $i\in I$.
  \item If $M=\bigcup_{i\in I}M_i$, then $\dim M=0$
    $\iff$ $\dim M_i=0$ for all $i\in I$.
  \item If $M=\colim_{i\in I} M_i$ is the colimit of a directed system,
    then
    \begin{equation*}
      \dim M\le \liminf_{i\in I} \dim M_i.
    \end{equation*}
  \end{enumerate}
\end{remark}

\begin{proof}[Proof of Theorem \ref{basic properties of capacity}]
\ref{submod}) Every
{\measurable} submodule of $M_0$ is a {\measurable} submodule of $M_1$, therefore
$\capac''(M_0)\le \capac''(M_1)$.

\ref{quotmod})
For every {\measurable} submodule $M$ of $M_2$ we find a finitely generated
submodule $N_M\subset M_1$ which projects onto $M$.
If $M_1$ is {\cofinalm},
then $N_M$ is {\measurable} and $\capac'(N_M)\ge \capac'(M)$ since on the left
we take the infimum over a smaller set of numbers. Therefore
$$\capac''(M_1)\ge \sup_{M\subset M_2\text{ {\measurable}}} \capac'(N_M)\ge
\sup_{M\subset M_2\text{ {\measurable}}}\capac'(M)=\capac''(M_2).$$

\ref{extension})
\newcounter{stepper}
Step \stepcounter{stepper}\arabic{stepper}: We prove $\capac(M_1)\le \capac(M_0)+\capac(M_2)$ if
$M_0$ is {\measurable} and $M_2$ is finitely presented. We will also see, that
this implies, that $M_1$ is {\measurable}:

By \cite[Lemma 3.4]{Lueck (1997a)} there is a finitely generated free
$\mathcala$-module
$F_2$  and an  exact sequence
\mbox{$0 \longrightarrow F_2 \stackrel{i_2}{\longrightarrow} F_2
  \longrightarrow M_2
\longrightarrow 0$}.
Let $q:F_0\to M_0$ be a projection of a finitely generated free
$\mathcala$-module onto $M_0$.

We get the following commutative diagram with exact rows and columns
where the $F_i$ are finitely generated free:
$$
\begin{CD}
& & 0 & & 0 & & 0 & &
\\
  & & @AAA @AAA @AAA
\\
0 @>>> M_0 @>>> M_1 @>>> M_2 @>>> 0
\\
  & & @AAA @AAA @AAA
\\
0 @>>> F_0 @>>> F_1 @>>> F_2 @>>> 0
\\
  & & @AAA @AAA @AAA
\\
0 @>>> K_0 @>>> K_1 @>>> i_2(F_2) @>>> 0
\\
  & & @AAA @AAA @AAA
\\
& & 0 & & 0 & & 0 & &
\end{CD}$$
Let  $K_0'\subset K_0$ be a finitely generated
submodule with $\dim F_0/K_0'=0$ (Lemma \ref{projinf}). We can
consider $K_0'$ also as submodule of $K_1$. Let
\mbox{$s: i_2(F_2) \longrightarrow K_1$}
be a section of the epimorphism
\mbox{$K_1 \longrightarrow i_2(F_2)$}.
Let $K_1'$ be the (finitely generated!) submodule of $K_1$
generated by $K_0'$ and the image of $s$. We obtain the exact
sequence
\begin{equation*}
  0\to K_0' \to K_1'\to i_2(F_2)\to 0 .
\end{equation*}
Going to the quotients, we obtain a commutative diagram with epimorphisms
as vertical maps whose lower row is an exact sequence
of finitely presented $\mathcala$-modules
$$\begin{CD}
0 @>>> M_0 @>>> M_1 @>>> M_2 @>>> 0
\\
  & & @AAA @AAA @AA\id A
\\
0 @>>> F_0/K_0' @>>> F_1/K_1' @>>> M_2 @>>> 0 .
\end{CD}$$
Note that this implies $M_1$ to be \measurable. By Lemma \ref{exact_seq_fp}
\begin{equation*}
  \capac(M_1)\stackrel{\ref{quotmod})}{\le}
  \capac(F_1/K_1')\le\capac(F_0/K_0')+\capac(M_2).
\end{equation*}
This holds for every $K_0'$ as above, therefore also
for the infimum $\capac(M_0)$ (by Lemma \ref{projinf})
in place of $\capac(F_0/K_0')$.

Step \stepcounter{stepper}\arabic{stepper}: We prove the inequality if $M_0$ and $M_1$ are arbitrary and
$M_2$ is finitely presented: Choose $N_1\subset M_1$ \measurable. Let
$N_0:=i^{-1}(N_1)$ and $N_2:=p(N_1)$. Let $L_1$ be a finitely presented
module with $\dim(L_1)=0$ projecting onto $N_1$. We get a commutative
diagram with exact
rows and surjective columns
\begin{equation*}
  \begin{CD}
    0 @>>> N_0 @>>> N_1 @>>> N_2 @>>> 0\\
    && @AAA   @AAA @AA\id A \\
    0 @>>> K_0 @>>> L_1 @>>> N_2 @>>>0.
  \end{CD}
\end{equation*}
$N_2$ is the image of the composition $L_1\to M_2$ and $K_0$ is the
kernel of this map. Since $L_1$ and $M_2$ are finitely presented, by
\cite[Theorem 0.2]{Lueck (1997a)} the same is true for $N_2$ and
$K_0$. In particular $N_0$ is \measurable. Therefore by Step 1
\begin{equation*}
  \capac(M_1)=\sup_{N_1\subset M_1 \text{ \measurable}}\capac(N_1)\le
  \sup(\capac(N_0)+\capac(N_2))\stackrel{\ref{submod})}{\le} \capac(M_0)+\capac(M_2).
\end{equation*}

Step \stepcounter{stepper}\arabic{stepper}: We prove the inequality if $M_1$ is \measurable: Obviously
$M_2$ is also measurable. Let $f:L_2\to M_2$ be a projection with
$L_2$ finitely presented and $\dim(L_2)=0$. By the pull back
construction we obtain a
commutative diagram with exact rows and
epimorphisms as vertical arrows
$$\begin{CD}
0 @>>> M_0 @>>> M_1 @>>> M_2 @>>> 0
\\
  & & @A \id AA @AAA @A f AA
\\
0 @>>> M_0 @>>> X_1 @>>> L_2 @>>> 0 .
\end{CD}$$
Note that $X_1$ as a submodule of $L_2\oplus M_1$ is {\cofinalm} by
the proof of Proposition \ref{capac well def}.\ref{capac' = capac''}.
Then by the last step and \ref{quotmod})
\begin{equation*}
  \capac(M_1)\stackrel{\ref{quotmod})}{\le} \capac(X_1)\le
  \capac(M_0)+\capac(L_2).
\end{equation*}
This holds for every $L_2$ as above, therefore also for the infimum $\capac(M_2)$.

Step \stepcounter{stepper} \arabic{stepper}: Finally $M_0$, $M_1$ and
$M_2$ are arbitrary $\mathcala$-modules. Suppose
$N_1\subset M_1$ is \measurable. We get the exact sequence $0\to
N_0\to N_1\to N_2\to 0$ as above with $N_i\subset M_i$. Then by the
previous step and \ref{submod})
\begin{equation*}
  \capac(N_1)\le
  \capac(N_0)+\capac(N_2)\stackrel{\ref{submod})}{\le}\capac(M_0)+\capac(M_2).
\end{equation*}
This holds for all $N_1$ as above and passing to the supremum yields
the desired inequality.

\ref{unions}.) This follows from \ref{submod}) and the fact,
  that each finitely generated and in particular each
{\measurable} submodule $L \subset M$ is contained in some $M_i$.

\ref{capac and colimits}.)
Let $N\subset M$ be a {\measurable} submodule. It suffices to show
\begin{equation*}
  \capac(N)\le \liminf_{i\in I} \capac(M_i).
\end{equation*}
Since $N$ is finitely generated and $M$ is the colimit of a directed
system, we find $i_0$ and a finitely generated $N_{i_0}\subset
M_{i_0}$ which projects onto $N$. For $i\ge i_0$ set
$N_i:=\phi_{i_0i}(N_{i_0})\subset M_i$. Then $\capac(N_i)\le
\capac(M_i)$ by \ref{submod}) and (since $I$ is directed) $\liminf_{i\ge
  i_0} \capac(N_i)\le \liminf_{i\in I}\capac(M_i)$. Observe that $N_i$
projects onto $N$ for all $i\ge i_0$. We will show that there is
$i_1\in I$
such that $N_i$ is {\measurable} for $i\ge i_1$. Then by \ref{quotmod})
$\capac(N)\le
\capac(N_i)$ for $i\ge i_1$ and therefore also $\capac(N)\le
\liminf_{i\ge i_1} \capac(N_i)$ which implies the  assertion.

Let $p_{i_0}:F\to N_{i_0}$
 be a projection  of a finitely generated free module onto $N_{i_0}$
 and let $p: F\to N$ be the composed projection.
By Lemma \ref{projinf}, since $N$ is {\measurable},
 we
find a finitely generated submodule $K\subset \ker(p)$  with $\dim(F/K)=0$.
 Because $K$
is finitely generated
and $I$ is a directed system, we find $i_1\ge i_0$ such that
$\phi_{i_0i_1}p_{i_0}(K)=0$. Then $\phi_{i_0i}p_{i_0}$ induces a
projection of the finitely presented zero-dimensional module
$F/K$ onto $N_i$ for every $i\ge i_1$. Therefore these $N_i$ are
{\measurable} and the  inequality follows.

For the second part suppose that every $M_i$ is {\measurable} and
every $\phi_{ij}$ is surjective. \ref{quotmod}) implies that $\liminf
\capac(M_i)=\inf \capac(M_i)$. It remains to prove $\capac(M)\ge\inf
\capac(M_i)$. Choose some $a\in I$ and a projection $p_a:F\to M_a$ with
$F$ finitely generated free. By composition we get a projection
$p:F\to M$. Let $K\subset \ker(p)$ be finitely generated with
$\dim(F/K)=0$. Since $I$ is a directed system, we find $b\ge a$ so
that  with $p_b:=\phi_{ab}p_a$ already $p_b(K)=0$. Since $\phi_{ab}$
is surjective, $p_b:F\to M_b$ is onto. Therefore by \ref{quotmod})
\begin{equation*}
  \capac(F/K)\ge \capac(M_b)\ge \inf\capac(M_i).
\end{equation*}
This holds for any finitely generated $K$ as above, therefore also for
the infimum in place of $\capac(F/K)$ which by Lemma \ref{projinf} is
$\capac(M)$.

\ref{capac and direct sums}.)
Since $\oplus_{i\in I}M_i = \bigcup_{\substack{J\subset I\\ \text{$J$
      finite}}} \oplus_{i\in J} M_i$  and because of \ref{unions}.) we
may assume $I$ is finite.
By induction we restrict to the case $\oplus_{i \in I} M_i=M_0 \oplus M_1$. Because of \ref{submod}) it remains to
prove
\begin{eqnarray}
\capac\left(M_0\oplus M_1\right) & \le &
\sup\{\capac(M_0),\capac(M_1)\}. \label{eqn 2.20}
\end{eqnarray}
If $M_0$ and $M_1$
are finitely presented \eqref{eqn 2.20} follows from \cite[Lemma
1.12]{Lott-Lueck (1995)} in the same way as Lemma \ref{exact_seq_fp}
does. If $M_0$ and $M_1$ are measurable,
we can choose epimorphisms $L_0 \to M_0$ and $L_1 \to M_1$ with $L_0$
and $L_1$ finitely presented and of dimension zero.
Note that $L_0 \oplus L_1$ is finitely presented and therefore
the result for the finitely presented
case implies $\capac(M_0 \oplus M_1) \leq \capac(L_0 \oplus L_1) \leq \sup \{ \capac(L_0), \capac(L_1) \}$. Since this
holds for
every choice of $L_0$ and $L_1$ we can pass to the infimum and get
\eqref{eqn 2.20}. Now let $M_0$ and $M_1$
be arbitrary modules. Then every measurable submodule $N \subset M_0
\oplus M_1 $ is contained in
$ N_0 \oplus N_1$ where $N_0$ and $N_1$ are the images  of $N$
under projection to $M_0$ and $M_1$. In particular they
are measurable as quotients of a measurable module.
Applying the result in the measurable case we get $\capac(N)=\capac(N_0 \oplus N_1) \leq \sup\{ \capac(N_0) ,
\capac(N_1) \}
\leq \sup \{ \capac(M_0) , \capac(M_1) \} $. Passing to the supremum on the left yields
\eqref{eqn 2.20} for arbitrary modules
$M_0$ and $M_1$. This finishes the proof of \ref{basic properties of capacity}.
\end{proof}

There are examples showing  that the  inequality
  \ref{extension})  in Theorem \ref{basic properties of capacity} is
  sharp. Moreover the assumption on cofinality in
\ref{quotmod}) is necessary by the following example.

\begin{example}
\label{finitely generated module with dim = 0 and kappa = 0^-}
We construct a non-trivial finitely generated $\mathcala$-module
$M$ with \mbox{$\dim(M) = 0$} which contains no
non-trivial {\measurable} $\mathcala$-submodule. In particular
$M$ is not {\cofinalm} and \mbox{$\capac(M) = 0^-$}.
Moreover, we construct
a quotient module $N$ of $M$ with
\mbox{$\capac(N) > \capac(M)$}.

Take \mbox{$\mathcala = L^{\infty}(S^1)$} which can be identified with
the group von Neumann algebra $\mathcaln(\zz)$. Let
$\chi_n$ be the characteristic function of the
subset $\{\exp(2\pi it)\mid t \in [1/n,1-1/n]\}$ of $S^1$.
Let \mbox{$P_n$} be the submodule in
$P=L^{\infty}(S^1)$ generated by $\chi_n$.
It is a direct summand. Hence the quotient $P/P_n$ is a finitely
generated  projective
$\mathcala$-module of dimension $2/n$. Projectivity implies
$\capac(P/P_n)=0^-$.
Define
\mbox{$I:=\bigcup_{n\in\nn} P_n\subset L^{\infty}(S^1)$}
and put \mbox{$M = L^{\infty}(S^1)/I$}.
Then $M=\colim_{n \in \nn } P/P_n$ is a finitely generated
$\mathcala$-module  with $\dim(M) = 0$ and $\capac(M)=0^-$ by
Theorem \ref{basic properties of capacity}.\ref{capac and
  colimits} and Remark \ref{basicdim}

Observe that the same argument applies to any von Neumann algebra $\mathcala$
with a directed system
$P_i\subset P$ of direct summands of a projective $\mathcala$-module
$P$ such that $\dim P_i<\dim P$ but $\dim
P=\sup_{i\in I} \dim P_i$.

For the quotient example,
put \mbox{$N_n = L^{\infty}(S^1)/((z-1)^n)$} for  positive integers
\mbox{$n \ge 0$}, where
$((z-1)^n)$ is the ideal generated by the function
\mbox{$S^1 \longrightarrow \cc$} sending $z$ to $(z-1)^n$.
Obviously \mbox{$I \subset ((z-1)^n)$} so that $N_n$ is a quotient of
$M$. The $\mathcala$-module $N_n$ is finitely presented with
\mbox{$\alpha(N_n) =1/n$} \cite[Example 4.3]{Lueck (1997a)}.
Hence we get
$$\capac(M) =  0^-;\quad\text{but}\quad \capac(N_n)  =  n. $$\qed
\end{example}

\begin{lemma}\label{elementary properties of special and
cofinal-special}
\begin{enumerate}
\item \label{L1}
A finitely generated submodule of a {\measurable} $\mathcala$-module
is again {\measurable}.

\item \label{L2}
A quotient module of a {\measurable} $\mathcala$-module is again {\measurable}.

\item \label{L4}
An $\mathcala$-module $M$ is {\cofinalm} if and only if
it is the union of its {\measurable} submodules.

\item \label{L5}
Submodules and  quotient modules of  {\cofinalm}
$\mathcala$-modules are again {\cofinalm}.

\item \label{L6}
Let \textshortexactsequence{M_0}{i}{M_1}{p}{M_2}
be an exact sequence of $\mathcala$-modules.
If $M_0$ and $M_2$ are {\cofinalm}, then $M_1$ is {\cofinalm}.

\item \label{L7}
The full subcategory of the abelian
category of all $\mathcala$-modules
consisting of {\cofinalm} modules is abelian and closed under
colimits over directed systems.
Given $r \in \{0^-\} \amalg [0,\infty]$, this is also
true for the full subcategory of {\cofinalm} $\mathcala$-modules $M$
with $\capac(M) \le r$.

\item \label{L8}
If $C_{\ast}$ is an $\mathcala$-chain complex
of {\cofinalm} $\mathcala$-modules,
then its homology \mbox{$H_p(C_{\ast})$}
is {\cofinalm} for all $p$. Moreover
$\capac(H_p(C)) \le \capac(C_p)$.

\end{enumerate}
\end{lemma}
\begin{proof}

\noindent
\ref{L1}.) This follows from the proof of Proposition \ref{capac well def}.

\noindent
\ref{L2}.) This is obvious.

\noindent
\ref{L4}.) Since $M$ is the  union of its finitely generated
submodules, it is the  union of its {\measurable} submodules
provided that $M$ is {\cofinalm}. Suppose that $M$
is the union of its {\measurable} submodules
and \mbox{$L\subset M$} is finitely generated.
There are finitely many
{\measurable} submodules $K_1, K_2, \ldots, K_r$
such that $L$ is contained in the submodule $K$
generated by $K_1, K_2, \ldots, K_r$. Obviously
$K$ and therefore $L$ is {\measurable} by assertions \ref{L1}.).

\noindent
\ref{L5}.) follows from assertionss \ref{L1}.) and  \ref{L2}.).

\noindent
\ref{L6}.) We have to show for a finitely generated submodule
\mbox{$M_1' \subset M_1$} that it is {\measurable}.
Let \mbox{$M_2' \subset M_2$} be the finitely generated submodule
\mbox{$p(M_1')$}  and \mbox{$M_0' \subset M_0$} be
\mbox{$i^{-1}(M_1')$}. Since $M_2$ is {\cofinalm},
$M_2'$ is {\measurable}. Choose a finitely presented
$\mathcala$-module $M_2''$ with \mbox{$\dim(M_2'') = 0$} together with an
epimorphism \mbox{$f: M_2'' \longrightarrow M_2'$}.
The pull back construction yields a commutative square
with exact rows and epimorphisms as vertical arrows
$$\begin{CD}
0 @>>> M_0' @>>> M_1' @>>> M_2' @>>> 0
\\
  & & @A \id AA @A\overline{f} AA @A f AA
\\
0 @>>> M_0' @>i''>> M_1'' @>p''>> M_2'' @>>> 0 .
\end{CD}$$
Since $M_1'$ is finitely generated, there is a finitely generated
submodule \mbox{$M_1''' \subset M_1''$}
such that \mbox{$\overline{f}(M_1''') = M_1'$}.
Let \mbox{$M_2''' \subset M_2''$} be the finitely generated submodule
\mbox{$p''(M_1''')$}  and \mbox{$M_0''' \subset M_0'$} be
\mbox{$(i'')^{-1}(M_1''')$}. We obtain an exact sequence
$$0 \longrightarrow M_0''' \longrightarrow M_1'''
\longrightarrow M_2''' \longrightarrow 0.$$
Since $M_2'''$ is a finitely generated
submodule of the finitely presented $\mathcala$-module
$M_2''$, the quotient \mbox{$M_2''/M_2'''$} is finitely presented.
Since the category of finitely presented $\mathcala$-modules is
abelian \cite[Theorem 0.2]{Lueck (1997a)},
the $\mathcala$-module $M_2'''$ is finitely presented.
Since $M_0'''$ is the kernel of an
epimorphism of the finitely generated
$\mathcala$-module $M_1'''$ onto
the finitely presented $\mathcala$-module $M_2'''$,
$M_0'''$ is finitely generated. As $M_0'''$ is a submodule of the
{\cofinalm} $\mathcala$-module $M_0$, the $\mathcala$-module $M_0'''$
is {\measurable}. Since $M_0'''$ is {\measurable} and
$M_2'''$ finitely presented, $M_1'''$ is {\measurable}
as follows from the first step of the proof of Theorem \ref{basic properties of
  capacity}.\ref{extension}). Hence $M_1'$ is {\measurable} by assertions \ref{L2}.),
since it is a quotient of $M_1'''$.

\noindent
\ref{L7}.) If $M=\colim_{i\in I}M_i$ is the colimit of a directed
system with $\psi_i:M_i\to M$, then $M$ is the directed union of
$\psi_i(M_i)$. If
all $M_i$ are {\cofinalm}, then their quotients $\psi_i(M_i)$ are
{\cofinalm} by assertions \ref{L5}.) and the same is true for $M$ by
assertions \ref{L4}.). The assertions now follows from \ref{L5}.),
\ref{L6}.) and Theorem
\ref{basic properties of capacity}.\ref{capacity and exact sequences}.

\noindent
\ref{L8}.) follows from \ref{L5}.) and Theorem
\ref{basic properties of capacity}.\ref{capacity and exact sequences}.
\end{proof}

Finally we discuss the behaviour of these notions under induction and
restriction for
subgroups $i:\Delta\to\Gamma$. The functor $i_*$ was
already investigated in \cite[Theorem 3.3]{Lueck (1998a)} where it is shown that $i_{\ast}$ is exact
and $\dim_{\mathcaln ( \Delta ) } (M)=\dim_{\mathcaln ( \Gamma )}(i_*M)$.

\begin{lemma} \label{capacity and induction}
Let $i: \Delta \longrightarrow \Gamma$ be an inclusion of groups,
then $i$ induces an inclusion $i:\mathcaln(\Delta)\to \mathcaln(\Gamma)$.
  \begin{enumerate}
\item \label{indgen}
If $M$ is a  {\measurable} or {\cofinalm}  respectively
$\mathcaln(\Delta)$-module,
then the $\mathcaln(\Gamma)$-module
$i_*M:=\mathcaln( \Gamma) \otimes_{\mathcaln(\Delta)} M$ is
{\measurable} or {\cofinalm}  respectively  and
\begin{eqnarray*}
\capac_{\mathcaln(\Delta)}(M) & = & \capac_{\mathcaln(\Gamma)}(i_*M).
\end{eqnarray*}
For an arbitrary $\mathcaln(\Delta)$-module $N$ we have
$\capac(N)\le\capac(i_*N)$.
\item \label{indfin} If the index of $\Delta$ in $\Gamma$ is finite
  and $M$ is an arbitrary $\mathcaln(\Delta)$-module, then
  \begin{equation*}
\capac_{\mathcaln(\Delta)}(M) = \capac_{\mathcaln(\Gamma)}(i_*M).
  \end{equation*}
\item \label{resfin}
If the index of $\Delta$ in $\Gamma$ is finite, $N$ is
 an $\mathcaln(\Gamma)$-module and $i^*N$ the
$\mathcaln(\Delta)$-module obtained by restriction, then
\begin{equation*}
  \capac_{\mathcaln(\Delta)}(i^*N)=\capac_{\mathcaln(\Gamma)}(N)
\end{equation*}
and $i^*N$ is {\measurable} or {\cofinalm} if
and only if $N$ has the same property.
\end{enumerate}
\end{lemma}
\proof
\noindent
\ref{indgen}.)
First suppose $M$ is a finitely presented zero-dimensional
$\mathcaln(\Delta)$-module.
Choose a resolution
$0\to F\stackrel{f}{\to} F\to M\to 0$ with a finitely generated free
module $F$ as in \cite[Lemma 3.4]{Lueck (1997a)}. Apply the proof  of
\cite[Lemma 3.6]{Lott-Lueck (1995)} to $f$, taking into account the
equivalence of free Hilbert $\mathcaln\Gamma$-modules and free
algebraic $\mathcaln\Gamma$-modules of \cite{Lueck (1997a)}. It
follows that $i_*M$ is a finitely presented $\mathcaln(\Gamma)$-module
with $\dim_{\mathcaln(\Gamma)}(i_*M) = 0$
and $\capac_{\mathcaln(\Delta)}(M)  = \capac_{\mathcaln(\Gamma)}(i_*M)$.

Next let $M$ be a {\measurable} $\mathcaln(\Delta)$-module.
Let $p:F \to M$ be a projection of a finitely generated free
$\mathcaln(\Delta)$-module onto $M$. Set $K:=\ker(p)$.
Then $i_*p:i_*F\to i_*M$ is surjective with kernel $i_*K$ (since $i_*$
is exact). If $K_1\subset
K$ is finitely generated with $\dim F/K_1=0$ (such a module exists by
Lemma \ref{projinf}),
then $i_*K_1\subset i_*K$ is also finitely generated with
\begin{equation*}
\dim
(i_*F/i_*K_1)\stackrel{\text{exactness}}{=} \dim i_*(F/K_1)=\dim F/K_1=0 .
\end{equation*}
Since $i_*F/i_*K_1$ is finitely presented and projects onto $i_*M$
the latter module is \measurable. Moreover by Lemma \ref{projinf} and
the first step applied to  $F/K_1$
\begin{equation*}
  \capac(i_*M)\le
  \inf_{K_1} \capac(i_*F/i_*K_1)=\inf_{K_1} \capac(F/K_1)=\capac(M) .
\end{equation*}

Choose on the other hand a finitely generated
$\mathcaln(\Gamma)$-submodule $L\subset i_*K$  with
$\dim (i_*F)/L=0$. For $L$ we find finitely many generators $\sum u_i \otimes k_i
\in i_*K = \mathcaln(\Gamma)\otimes_{\mathcaln(\Delta)} K$. Let
$L'$ be the submodule of $K$ generated by all the $k_i$. Then $L\subset
i_*L'$, therefore
\begin{equation*}
  0\le \dim F/L' = \dim i_*F/i_*L' \le \dim i_*F/L =0
\end{equation*}
and (since $F/L'$ is finitely presented and by
\ref{basic properties of capacity}.\ref{quotmod}))
\begin{equation*}
  \capac(M)\le \capac(F/L')=\capac(i_*F/i_*L')\le \capac(i_*F/L) .
\end{equation*}
Since this holds for arbitrary $L$ as above, Lemma \ref{projinf}
implies
$\capac(M)\le \capac(i_*M)$.

If $M$ is {\cofinalm}, then it is the union $\cup_{i \in I} M_i$ over
the directed
system of its {\measurable} $\mathcaln(\Delta)$-submodules. Since $i_*$ is exact,
the $\mathcaln(\Gamma)$-module $i_*M$ is the union
$\bigcup_{i \in I} i_*M_i$ over the directed
system of {\measurable} $\mathcaln(\Gamma)$-submodules $i_*M_i$.
We conclude from Theorem
\ref{basic properties of capacity}.\ref{unions} and the previous step
that $i_*M$ is {\cofinalm} and
\begin{equation*}
\capac_{\mathcaln(\Delta)}(M) =
\sup_{i\in I}\{\capac_{\mathcaln(\Delta)}(M_i)\}
=
\sup_{i\in
  I}\{\capac_{\mathcaln(\Gamma)}(i_*M_i))\}
=\capac_{\mathcaln(\Gamma)}(i_*M) .
\end{equation*}

Last let $M$ be an arbitrary $\mathcaln(\Delta)$-module. Since every
{\measurable} $\mathcaln(\Delta)$-submodule of $M$ induces
an $\mathcaln(\Gamma)$-submodule of $i_*M$ of the same capacity,
$\capac(M)\le\capac(i_*M)$ by Definition \ref{definition of capacity}.

\noindent\ref{resfin}.)
We begin with studying the restriction. Here $i^*\mathcaln(\Gamma) =
\oplus_{i=1}^{[\Gamma:\Delta]} \mathcaln(\Delta)$ since the same holds
for $\cc\Gamma$ as a $\cc\Delta$-module and
$\mathcaln(\Gamma)=\mathcaln(\Delta)\tensor_{\cc\Delta}\cc\Gamma$.
This observation and the proof of \cite[Lemma 3.6]{Lott-Lueck
  (1995)} imply that if $N$ is a finitely presented
$\mathcaln(\Gamma)$-module, then $i^*N$ is finitely presented as
$\mathcaln(\Delta)$-module and
\begin{equation*}
  \dim_{\mathcaln(\Delta)}(i^*N)=[\Gamma:\Delta]\dim_{\mathcaln(\Gamma)}(N);\qquad \capac_{\mathcaln(\Delta)}(i^*N) =
\capac_{\mathcaln(\Gamma)}(N).
\end{equation*}
If $N$ is arbitrary and $L\to N$ a projection of a finitely presented
zero-dimensional $\mathcaln(\Gamma)$-module, then $i^*L\to i^*N$ is a
corresponding projection and $\capac(L)=\capac(i^*L)$. If on the other
hand $L'$ is a zero-dimensional $\mathcaln(\Delta)$-module projecting
onto $i^*N$, then
$i_*L'=\mathcaln(\Gamma)\tensor_{\mathcaln(\Delta)}L'$ is a finitely
presented $\mathcaln(\Gamma)$-module
naturally projecting onto $N$ with the same  dimension and capacity. In
particular $N$ is {\measurable} if and only if $i^*N$ is {\measurable}
and by Definition \ref{definition of capacity} the capacities coincide
in this case.

Any {\measurable} submodule of an $\mathcaln(\Gamma)$-module $N$ restricts
to a {\measurable} $\mathcaln(\Delta)$-submodule of $i^*N$ with the same
capacity. On the other hand, if $U\subset i^*N$ is a {\measurable}
$\mathcaln(\Delta)$-submodule and $V$ is the
$\mathcaln(\Gamma)$-module generated by $U$, then $V$ is a quotient of
$i_*U$ and U a submodule of V, therefore by assertions \ref{indgen}.) and
Theorem
\ref{basic properties of capacity}.\ref{capacity and exact sequences}
$V$ is {\measurable} and
\begin{equation*}
\capac(U)=\capac(i_*U)\ge
\capac(V)=\capac(i^*V) \ge \capac(U).
\end{equation*}
Definition \ref{definition of capacity} implies
$\capac(N)=\capac(i^*N)$. The {\cofinalm} case is proven as above.

\noindent\ref{indfin}.)
Let $M$ be an $\mathcaln(\Delta)$-module. Then $i^*i_*M \cong
\oplus_{i=1}^{[\Gamma:\Delta]} M$. By \ref{indgen}.) and
\ref{resfin}.)
\begin{equation*}
  \capac(M)\le \capac(i_*M) = \capac(i^*i_*M) =
  \capac(\oplus_{i=1}^{[\Gamma:\Delta]}M)=\capac(M).
\end{equation*}\qed

\typeout{--------------------- section 2 -----------------------}
\section{Capacity of groups}
\label{Capacity of groups}

In this section, we apply the concepts developed so far to define and
study Novikov-Shubin invariants respectively capacities to arbitrary
spaces, via the classifying space in particular to arbitrary groups.

\begin{definition} \label{definition of capacity for spaces and groups}
Let $X$ be a topological space with an action of the discrete group $\Gamma$.
Let $H^{\Gamma}_p(X;\mathcaln(\Gamma))$ be the $\mathcaln(\Gamma)$-module given by
the $p$-th homology of the chain complex
$\mathcaln(\Gamma) \otimes_{\zz\Gamma} C^{\sing}_*(X)$,
where $C^{\sing}_*(X)$ is the singular chain complex of $X$.
Define the {\em  $p$-th capacity} of $X$
$$\capac_p(X;\mathcaln(\Gamma)) ~ := \capac(H^{\Gamma}_p(X;\mathcaln(\Gamma))).$$
Define the {\em  $p$-th capacity} of the group $\Gamma$ by
$$\capac_p(\Gamma) ~ := ~ \capac_p(E\Gamma;\mathcaln(\Gamma)), $$
where $E\Gamma$ is any universal free $\Gamma$-space.
\end{definition}

A group is {\em locally } finite, nilpotent, abelian, $\ldots$
respectively, if any finitely generated subgroup is
finite, nilpotent, abelian, $\ldots$ respectively.
A group is {\em virtually} nilpotent, abelian, $\ldots$
respectively, if it contains a subgroup of finite index which is
nilpotent, abelian, $\ldots$ respectively.

If $S$ is a finite set of generators for the group $G$,
let $b_S(k)$ be the number of elements in $G$ which can be written
as a word in $k$ letters of $S \cup S^{-1} \cup \{1\}$. The group $G$ has
\emph{polynomial growth of degree not greater than $n$}
if there is $C$ with $b_S(k) \le C \cdot k^d$ for all $k \ge 1$.
This property is a property of $G$ and independent of the choice of the
finite set $S$ of generators.
We say that $G$ has \emph{polynomial growth}
if it has polynomial growth of degree not greater than $n$ for some $n > 0$.
A finitely generated
group $\Gamma$ is \emph{nilpotent}
if $\Gamma$ possesses a finite \emph{lower central series}
$$\Gamma = \Gamma_1 \supset \Gamma_2 \supset \ldots \supset \Gamma_s = \{1\}
\hspace{15mm} \Gamma_{k+1} = [\Gamma,\Gamma_{k}].$$
Let $n_i$ be the rank of the finitely generated abelian group
$\Gamma_i/\Gamma_{i+1}$ and let $n$
be the integer $\sum_{i\ge 1} i \cdot n_i$. Suppose that $\overline{\Gamma}$
contains $\Gamma$ as subgroup of finite index.
Then for any
finite set $S$ of generators of $\overline{\Gamma}$
there is a constant $C > 0$
such that  $C^{-1} \cdot k^{n} \le b_S(k) \le C\cdot k^{n}$ holds
for any $k \ge 1$ and in particular
$\overline{\Gamma}$ has
polynomial growth precisely of degree $n$
\cite[page 607 and Theorem 2 on page 608]{Bass(1972)}.

A famous result of Gromov \cite{Gromov (1981)}
says that a finitely generated group has
polynomial growth if and only if it is virtually nilpotent.

\begin{proposition} \label{zero-th capacity}
\begin{enumerate}

\item \label{H_0(Gamma) cofinal-special}\label{c0fgappr}
If $\Gamma$ is not locally finite, then
$H^{\Gamma}_0(E\Gamma;\mathcaln(\Gamma))$ is {\measurable} and
\begin{equation*}
\capac_0(\Gamma)= \inf \{ \capac_0(\Gamma')|\; \Gamma'\subgroup\Gamma\text{
     infinite finitely generated} \} .
\end{equation*}

If $\Gamma$ is locally finite, then $\capac_0(\Gamma)= 0^-$
but $H^{\Gamma}_0(E\Gamma;\mathcaln(\Gamma))$
is non-trivial and is in particular not
\cofinalm.

\item \label{capac for f.g. v. nilpotent}
Suppose that $\Gamma$ is finitely generated. Then
\begin{equation*}
  \capac_0(\Gamma) =
  \begin{cases}
    0^-; &\text{if $\Gamma$ is finite or non-amenable},\\
    1/n; & \text{if $\Gamma$ has
      polynomial growth of
      degree $n$,}\\
    0; & \text{if $\Gamma$ is infinite and amenable, but not virtually
      nilpotent}.
  \end{cases}
\end{equation*}
  This
  computes $\capac_0$ for every finitely generated group.

\item \label{capac = 0^-}\label{capac not= 0^-,0}
Let $\Gamma$ be an arbitrary group. Then
\begin{equation*}
 \text{$\Gamma$ is locally finite or non-amenable} \iff \capac_0(\Gamma)=
  0^- .
\end{equation*}
If $\Gamma$ is locally
        virtually nilpotent but not locally finite, then
        \begin{equation*}
    \capac_0(\Gamma) = \inf\{\capac_0(\Gamma')|\; \Gamma'\subgroup\Gamma\text{ infinite finitely generated
        nilpotent} \} .
  \end{equation*}
If $\Gamma$ is amenable and contains a subgroup which is finitely generated but
        not virtually nilpotent, then
\begin{equation*} \capac_0 (\Gamma)= 0.
\end{equation*}
Note that every group belongs to
        one of the categories and that $c_0(\Gamma)>0$ implies that
        $\Gamma$ is locally virtually nilpotent but not locally finite.

\end{enumerate}
\end{proposition}
The above also applies to $c_0$ of arbitrary
path-connected $\Gamma$-spaces by \cite[4.10]{Lueck (1998a)}.
Observe that $H_0^{\Gamma}( E \Gamma ; \mathcaln ( \Gamma ))$  is not
\cofinalm\ if $\Gamma$ is finite or locally finite. This is
responsible for the clumsiness of some of the statements below
because Theorem \ref{basic properties of capacity}
\ref{quotmod}) becomes false without the condition \cofinalm\
as shown in Example
\ref{finitely generated module with dim = 0 and kappa = 0^-}.

\begin{proof}
Remember that $H^\Gamma_0(E\Gamma;\mathcaln(\Gamma)) =
\mathcaln(\Gamma)\tensor_{\cc\Gamma}\cc$
which has dimension zero if and only if $\Gamma$ is infinite and is
trivial if and
only if $\Gamma$ is nonamenable by \cite[4.10]{Lueck
  (1998a)}. Moreover $\mathcaln(\Gamma)\tensor_{\cc\Gamma}\cc$ is finitely presented if
$\Gamma$ is finitely generated.

\noindent
\ref{c0fgappr}.)
If $\Gamma^{\prime}$ is finitely generated infinite
$\mathcaln(\Gamma^{\prime})\tensor_{\cc\Gamma^{\prime}}\cc$
is measurable and by Lemma
\ref{capacity and induction}.\ref{indgen}
$\mathcaln(\Gamma) \otimes_{\mathcaln( \Gamma^{\prime})} \mathcaln(\Gamma^{\prime})\tensor_{\cc\Gamma^{\prime}}\cc$ is
measurable and
$c(\mathcaln(\Gamma) \otimes_{\mathcaln( \Gamma^{\prime})} \mathcaln(\Gamma^{\prime})\tensor_{\cc\Gamma^{\prime}}\cc)=
c(\mathcaln(\Gamma^{\prime})\tensor_{\cc\Gamma^{\prime}}\cc)$. If $\Gamma$ is not locally
finite the system of infinite finitely generated subgroups is cofinal and therefore
\begin{equation*}
  \mathcaln(\Gamma)\tensor_{\cc\Gamma}\cc = \colim
  \mathcaln(\Gamma)\tensor_{\cc\Gamma'}\cc ,
\end{equation*}
where the colimit is taken over the directed system of infinite finitely
generated subgroups. Now the second part of Theorem
\ref{basic properties of capacity}.\ref{capac and colimits} yields the claim
if $\Gamma$ is not locally finite.

The proof of Theorem
\ref{estimates on capac, extensions and local information}.%
\ref{capac and finitely generated subgroups} shows for locally finite $\Gamma$
$$H^{\Gamma}_0(E\Gamma;\mathcaln(\Gamma)) ~ = ~
\colim_{\Delta \subset \Gamma} H^{\Delta}_0(E\Delta;\mathcaln(\Gamma))
~ = ~
\colim_{\Delta \subset \Gamma} H_0(E\Delta) \otimes_{\cc \Gamma}
\mathcaln(\Gamma)),$$
where $\Delta \subset \Gamma$ runs though the finite subgroups.
For finite $\Delta$, the $\cc \Delta$-module $H_0(E\Delta)$ is projective
and hence the $\mathcaln(\Gamma)$-module
$H^{\Delta}_0(E\Delta;\mathcaln(\Gamma))$ is projective which implies
$\capac(H^{\Delta}_0(E\Delta;\mathcaln(\Gamma))) = 0^-$. By
Theorem \ref{basic properties of capacity}.\ref{capac and colimits}
$$\capac(\Gamma) = \capac(H^{\Gamma}_0(E\Gamma;\mathcaln(\Gamma))) = 0^-.$$
Since $\Gamma$ is non-amenable, $H^{\Gamma}_0(E\Gamma;\mathcaln(\Gamma))$
is non-trivial and hence cannot be \cofinalm.
\noindent

\ref{capac for f.g. v. nilpotent}.)
If $\Gamma$ is finite
or non-amenable, then
 $\mathcaln(\Gamma)\tensor_{\cc\Gamma}\cc$ is finitely generated projective
or trivial and therefore $\capac_0(\Gamma)=0^-$. If $\Gamma$ is infinite and amenable,
then
$H_0^\Gamma(E\Gamma;\mathcaln(\Gamma))$ is finitely presented and zero-dimensional but
non-trivial, therefore $\capac_0(\Gamma)\ge 0$.
The rest is the content of Lemma \ref{c0varop}.

\noindent
\ref{capac = 0^-}.)
This follows by combining the above results.
\end{proof}

\begin{lemma}\label{c0varop}
  Suppose $\Gamma$ is a finitely generated infinite group. Then
  $\capac_0(\Gamma)=1/n$ if $\Gamma$ has polynomial growth precisely of
  rate $n$. If $\Gamma$ is not virtually nilpotent, then
  $\capac_0(\Gamma) \le 0$.
\end{lemma}
\begin{proof}
  First observe that $\Gamma$ has polynomial growth of precisely degree
$n$ if and only if  the recurrency probablity $p(k)$ of
the natural
random walk on $\Gamma$ decreases polynomially with exponent $n/2$,
i.e. there is a constant $C > 0$ with
$C^{-1}\cdot k^{-n/2} \le p(n) \le C\cdot k^{-n/2}$ for $k \ge 1$
and it is not virtually nilpotent if and only if for any
$n > 0$ there is $C(n) > 0$ satisfying
$p(k) \le C \cdot k^{-n}$ for $k \ge 1$.
These are results of Varopoulos
\cite{Varopoulos(1988)}, compare \cite[6.6 and 6.7]{Woess (1994)}.

Now we have to translate this statement
to information about the spectrum of the Laplacian.
We will prove the following result from which our lemma follows:
  the finitely generatef group $\Gamma$ has
  polynomially decreasing recurrence probability with exponent
  $n/2$ if and only if $\capac_0(\Gamma)=1/n$.

  This can be deduced from
  \cite{Varopoulos(1988)}. We will give a self-contained and simple
  proof.

  If $S$ is a finite set of generators of $\Gamma$,  then
  $$H_0^{\Gamma}(E\Gamma;\mathcaln(\Gamma))= \coker(
  \oplus_{s\in S} \mathcaln(\Gamma)\xrightarrow{d_0=\oplus
    (s-1)}\mathcaln(\Gamma)) . $$
  Therefore $\capac_0(\Gamma)$ is the inverse of the Novikov-Shubin
  invariant $\alpha_0(\Gamma)$ of $d_0$, which we compute from
  the combinatorial Laplacian $\Delta_0=1-P$. Here $P$  is the
  transition operator $P(g)=(1/|S|)\cdot \sum_{s\in
    S}sg$. The recurrence probability is given
  by
  \begin{equation*}
    p(k):= (P^k(e),e) = \tr(P^k),
  \end{equation*}
and the spectral density function of $\Delta_0$ by
\begin{equation*}
  F(\lambda)=\tr(\chi_{[1-\lambda,1]}(P)).
\end{equation*}
 All of the operators in question are positive and
  therefore for $k\in\nn$ and $0<\lambda<1$
  \begin{equation}\label{opieq}
   (1-\lambda)^k\chi_{[1-\lambda,1]}(P)\le P^k\le
   (1-\lambda)^k\chi_{[0,1-\lambda]}(P) + \chi_{[1-\lambda,1]}(P) \le 1.
 \end{equation}
 Application of the trace to these inequalities gives
\begin{equation}
   \label{trieq}
   (1-\lambda)^k F(\lambda)\le p(k)\le (1-\lambda)^k + F(\lambda)
   \qquad \mbox{ for all } 0<\lambda<1.
 \end{equation}
 The first inequality implies if $0<\lambda<1$
 \begin{equation*}
   \frac{\ln F(\lambda)}{\ln\lambda}\ge \frac{\ln
     p(k)}{\ln\lambda}-k\frac{\ln(1-\lambda)}{\ln\lambda}.
 \end{equation*}
 If $p(k)\le Ck^{-a}$ for $C>0$, then, putting $k=[\lambda^{-1}]$ (the
 largest integer not larger than $\lambda^{-1}$) we see
 \begin{equation*}
   \alpha_1(\Gamma) =  2\liminf_{\lambda\to0^+}\frac{\ln
     F(\lambda)}{\ln\lambda}\ge
   2\lim_{\lambda\to 0}\left( a
   \frac{\ln\lambda}{\ln\lambda} +\frac{\ln C}{\ln\lambda}
   -\frac{\ln(1-\lambda)}{\lambda\ln\lambda}\right)=  2a .
\end{equation*}

Suppose now that $p(k)\ge C k^{-a}$. Choose $\epsilon>0$. Putting
$k:=[\lambda^{-(1+\epsilon)}]+1$ the second part of \eqref{trieq} implies
\begin{equation*}
  F(\lambda)\ge C\lambda^{a(1+\epsilon)}
  \left(\frac{[\lambda^{-1-\epsilon}]+1}{\lambda^{-1-\epsilon}}\right)^{-a}
    -(1-\lambda)^{[\lambda^{-(1+\epsilon)}]+1}.
\end{equation*}
Using lemma \ref{ie} below with $\delta=C/2$ and $a(1+\epsilon)$ instead of
$a$ this implies
\begin{equation*}
  \begin{split}
    &(1-\lambda)^{[\lambda^{-1-\epsilon}]+1}\stackrel{1-\lambda<1}{\le}(1-\lambda)^{-1-\epsilon}\le
\frac{C}{2}\lambda^{a(1+\epsilon)}\\
\implies &
    \frac{\ln F(\lambda)}{\ln\lambda}\le \underbrace{\frac{\ln(C
        \left(([\lambda^{-1-\epsilon}]+1)/\lambda^{-1-\epsilon}\right)^{-a}-C/2)}{\ln
        \lambda}}_{\text{$\to 0$ if $\lambda\to 0$}}
    + a(1+\epsilon)\frac{\ln\lambda}{\ln\lambda}.
\end{split}
\end{equation*}

Since the inequality is true for arbitrary $\epsilon>0$ we conclude
\begin{equation*}
  \alpha_0(\Gamma) = 2\liminf_{\lambda\to 0^+}\frac{\ln
    F(\lambda)}{\ln\lambda}\le 2a.
\end{equation*}
Now Lemma \ref{c0varop} follows.
\end{proof}

\begin{lemma}\label{ie}
  For arbitrary $\epsilon,\delta,a>0$ one finds $\lambda_0>0$ so that
  \begin{equation*}
    (1-\lambda)^{\lambda^{-1-\epsilon}}\le
    \delta\lambda^a\qquad\mbox{ for all } 0<\lambda<\lambda_0.
  \end{equation*}
\end{lemma}
\begin{proof}
  Note that for $0<\lambda<1$ the stated inequality is equivalent to
  \begin{equation*}
    \begin{split}
      &\lambda^{-1-\epsilon}\ln(1-\lambda) \le \ln\delta +
      a\ln\lambda\\
      \iff  & 1 \ge \frac{(\ln\delta +a
        \ln\lambda)\lambda^{1+\epsilon}}{\ln(1-\lambda)}.
    \end{split}
    \end{equation*}
    For $\lambda\to 0$ the right hand side tends to $0$ which can be
    seen using l'Hospital's rule.
\end{proof}

\begin{theorem} \label{estimates on capac, extensions and local information}
\begin{enumerate}

\item \label{capac and subgroups of finite index}
Let $\Delta \subset \Gamma$ be a subgroup of finite index.
Then
\begin{eqnarray*}
\capac_n(\Gamma) & = & \capac_n(\Delta).
\end{eqnarray*}
Moreover $H^\Delta_n(E\Delta;\mathcaln(\Delta))$ is {\cofinalm} if and only if
$H^\Gamma_n(E\Gamma;\mathcaln(\Gamma))$ is {\cofinalm}.

\item
\label{capac and normal subgroups}
Let $\Delta \subset \Gamma$ be a normal subgroup. Suppose that
$H^{\Delta}_q(E\Delta;\mathcaln(\Delta))$ is {\cofinalm} for $q \le n$.
Then we get for $p = 0,1,  \ldots, n$
that $H^{\Gamma}_p(E\Gamma;\mathcaln(\Gamma))$ is {\cofinalm} and
\begin{eqnarray*}
\capac_p(\Gamma) & \le & \sum_{q = 0}^p \capac_q(\Delta).
\end{eqnarray*}

\item \label{capac and finitely generated subgroups}
If there is a cofinal system of subgroups $\Delta \subset
\Gamma$ with
$H_p^{\Delta}(E \Delta ; \mathcaln ( \Delta))$
\cofinalm, then $H_p^{\Gamma}(E \Gamma ; \mathcaln ( \Gamma ))$
is \cofinalm and
\begin{eqnarray*}
c_n(\Gamma) \le  \liminf \{ c_n(\Delta) \},
\end{eqnarray*}
where $\Delta$ runs over the cofinal system.

\item
  \label{Z^n}
  If $n\ge 1$ and  $\Gamma=\zz^n$ or $\zz^{\infty}=\oplus_{i=1}^{\infty} \zz$,
   then
  $H_p^\Gamma(E\Gamma;\mathcaln(\Gamma))$ is {\cofinalm} for every $p$
 and
  \begin{equation*}
    \begin{split}
    \capac_p(\zz^n) & =
    \begin{cases}
      1/n & \mbox{ if }\; 0\le p\le n-1; \\ 0^- & \mbox{ if } \; p\ge n;
    \end{cases} \\
    \capac_p(\zz^\infty)  & \le 0 \quad \mbox{ if }  \; p\ge 0.
  \end{split}
  \end{equation*}
 (Remark: It is possible to show $\capac_p(\zz^\infty)=0$ for all
  $p\ge 0$.)
 \item \label{capac_1 for virtually nilpotent groups}
Suppose that $\Gamma$ is a finitely generated
virtually nilpotent group but not finite. Then
$H^{\Gamma}_p(E\Gamma;\mathcaln(\Gamma))$ is {\cofinalm} for $p \ge 0$ and
\begin{eqnarray*}
\capac_0(\Gamma) + \capac_1(\Gamma) & \le & 1;
\\
\capac_p(\Gamma) & \le & 1 \quad \text{ for } p \ge 1.
\end{eqnarray*}
The same holds for $\Gamma$ locally virtually nilpotent but not locally finite.
\end{enumerate}
\end{theorem}
\proof
\ref{capac and subgroups of finite index}.)
Let $i: \mathcaln(\Delta) \longrightarrow \mathcaln(\Gamma)$
be the ring homomorphism
induced by the inclusion. Since
$\mathcaln(\Delta) \otimes_{\zz\Delta} \zz\Gamma$ is isomorphic to  $\mathcaln(\Gamma)$
as
$\mathcaln(\Delta)$-$\zz\Gamma$-bimodule  and since $E\Gamma$ viewed as
$\Delta$-space is a model for $E\Delta$, we get
$i^*H^{\Gamma}_n(E\Gamma;\mathcaln(\Gamma)) = H^{\Delta}_n(E\Delta;\mathcaln(\Delta))$.
The statement now follows from Lemma \ref{capacity and induction}.

\noindent
\ref{capac and normal subgroups}.) There is a spectral sequence converging
to $H^{\Gamma}_{p + q}(E\Gamma;\mathcaln(\Gamma))$
whose $E_1$-term is given by
$$E_{p,q}^1 =  H^{\Delta}_q(E\Delta;\mathcaln(\Gamma)) \otimes_{\zz\pi}C_p(E\pi) =
\oplus_{I_p} i_*H^{\Delta}_q(E\Delta;\mathcaln(\Delta)),$$
where $i: \Delta \longrightarrow \Gamma$ is the inclusion and
$I_p$ the set of $p$-cells of $B\pi$. We conclude from
Theorem \ref{basic properties of capacity}.\ref{capac and direct sums}
and Lemma \ref{capacity and induction}.\ref{indgen}
that $E^1_{p,q}$ is {\cofinalm} for $q \le n$ and
$\capac(E^1_{p,q})  =  \capac_q(\Delta)$.
We conclude from Lemma \ref{elementary properties of special and
cofinal-special}.\ref{L8}
that $E^{\infty}_{p,q}$ is {\cofinalm} for $q \le n$ and
$\capac(E^{\infty}_{p,q}) \le   \capac_q(\Delta)$.
Theorem \ref{basic properties of capacity}.\ref{capacity and exact sequences}
and \ref{elementary properties of special and cofinal-special}.\ref{L6}
implies that $H^{\Gamma}_q(E\Gamma;\mathcaln(\Gamma))$ is {\cofinalm} for
$q \le n$ and
$\capac_q(\Gamma) \le \sum_{p=0}^q \capac_p(\Delta)$
for $0 \le p \le q$.

\noindent
\ref{capac and finitely generated subgroups}.)
Since $\Gamma$ is the union of the subgroups $\Delta$,
one can choose a model for $E\Gamma$ such that for each
subgroup $\Delta$ the model $E\Delta$ is a
subcomplex of $E\Gamma$ and $E\Gamma$ is the union of all the $E\Delta$'s.
For instance take as model for $E\Gamma$ the infinite join
$\Gamma \ast \Gamma \ast \dots$. Hence
\begin{eqnarray*}
E\Gamma & = &
\colim_{\Delta\subset\Gamma} \Gamma\times_{\Delta} E\Delta;
\\
H^{\Gamma}_p(E\Gamma;\mathcaln(\Gamma)) & =  &
\colim_{\Delta \subset \Gamma} H^{\Delta}_p(E\Delta;\mathcaln(\Gamma)),
\end{eqnarray*}
where $\Delta \subset\Gamma$ runs through the finitely generated subgroups.
Since $H^{\Delta}_p(E\Delta;\mathcaln(\Gamma))$ is $\mathcaln(\Gamma)$-isomorphic to
$\mathcaln(\Gamma) \otimes_{\mathcaln(\Delta)} H^{\Delta}_p(E\Delta;\mathcaln(\Delta))$
the claim follows from Theorem
\ref{basic properties of capacity}.\ref{capac and colimits} and
Lemma \ref{capacity and induction}.\ref{indgen}.

\noindent
\ref{Z^n}.)
  A direct computation shows the result for $\Gamma=\zz^n$. For $\zz^\infty$ apply assertions \ref{capac and finitely
  generated subgroups}.

\noindent \ref{capac_1 for virtually nilpotent groups}.) By \ref{capac and finitely generated subgroups}.)
we can assume that $\Gamma$ is finitely generated infinite and virtually
nilpotent. We claim that such a group $\Gamma$
is either virtually abelian
or contains a normal subgroup $\Delta$ such that
there exists a central extension
$1 \longrightarrow \zz \longrightarrow
\Delta \longrightarrow \zz^2
\longrightarrow 1$.
This is proven as follows.

Recall that subgroups and quotient groups of
nilpotent groups are nilpotent again,
any nilpotent group contains a normal torsionfree group of finite index
and the center of a non-trivial nilpotent group is non-trivial.
Now choose a normal torsionfree subgroup $\Gamma_0$ of $\Gamma$
of finite index and inspect the exact sequence
$1 \longrightarrow \cent(\Gamma_0) \longrightarrow
\Gamma_0 \longrightarrow \Gamma_0/\cent(\Gamma_0)
\longrightarrow 1$. If $\Gamma_0/\cent(\Gamma_0)$ is finite,
$\Gamma$ is virtually abelian. Suppose that
$\Gamma_0/\cent(\Gamma_0)$ is infinite. By inspecting the analogous
sequence for a normal torsionfree subgroup
$\Gamma_1 \subset \Gamma_0/\cent(\Gamma_0)$ of finite index and using
the fact that $\Gamma_1/\cent(\Gamma_1)$ is either finite or contains
$\zz$ as normal subgroup, one sees that
$\Gamma_0/\cent(\Gamma_0)$ contains $\zz$ as subgroup of
finite index or contains both $\zz$ and $\zz^2$ as normal subgroups.
Since $\cent(\Gamma_0)$ has at least rank $1$, the claim follows
for $\Gamma_0$ and hence for $\Gamma$.

If $\Gamma$ is virtually abelian
the assertion \ref{capac_1 for virtually nilpotent groups}.) follows from
\ref{capac and subgroups of finite index}.) and
\ref{Z^n}.). Suppose that $\Delta$ is a normal subgroup of $\Gamma$
and that there is a
central extension
$1 \longrightarrow \zz \longrightarrow
\Delta \longrightarrow \zz^2
\longrightarrow 1$.
Because of \ref{capac and normal subgroups}.)
$\capac_0(\Gamma) + \capac_1(\Gamma) ~ \le ~
2 \cdot \capac_0(\Delta) + \capac_1(\Delta).$
Hence it remains to show
$$2 \cdot \capac_0(\Delta) + \capac_1(\Delta) \le 1.$$
One can realize $\Delta$ as the fundamental group of a closed
$3$-manifold $M$ which is a principal $S^1$-bundle over $T^2$.
Hence $M$ is a Seifert manifold whose base orbifold has
Euler characteristic $\chi=0$ and the computation in
\cite[Theorem 4.1]{Lott-Lueck (1995)} shows
$\capac_0(\Delta) = \capac_1(\Delta) = 1/3$ if the Euler class $e(M)  = 0$
and $\capac_0(\Delta) = 1/4$ and $\capac_1(\Delta) = 1/2$
if $e(M) \not=0$. This finishes the proof of
Theorem \ref{estimates on capac, extensions and local information}.
\qed\par

\begin{definition}
  Given a class of groups ${\mathcal{X}}$, let $L{\mathcal X}$ be the
  class of groups $\Gamma$ for which any finitely generated subgroup
  $\Delta$ belongs to ${\mathcal X}$.  Given classes of groups
  ${\mathcal X}$ and ${\mathcal Y}$, let ${\mathcal X}{\mathcal Y}$ be
  the class of groups $\Gamma$ which contain a normal subgroup $\Delta
  \subset \Gamma$ with $\Delta \in {\mathcal X}$ and $\Gamma/\Delta
  \in {\mathcal Y}$. The {\em class ${\mathcal E}$ of elementary
    amenable groups} is defined as the smallest class of groups which
  contains all abelian and all finite groups, is closed under
  extensions, taking subgroups, forming quotient groups and under
  directed unions. The class of finite groups is denoted by $\mathcal F$.
\end{definition}

\begin{theorem} \label{estimate on capac_1 of an elementary amenable group}
Let ${\mathcal C}$ be the class of groups $\Gamma$ for which
$H^{\Gamma}_p(E\Gamma ;\mathcaln(\Gamma))$ is {\cofinalm} for all $p \ge 0$ and
$\capac_0(\Gamma) + \capac_1(\Gamma) \le 1$ holds.

\begin{enumerate}
\item \label{not in cal C}
Finite and locally finite groups do {\bf not} belong to $\mathcalc$.

\item \label{cal C and cal E}
If the infinite elementary amenable group $\Gamma$ contains no
infinite locally finite subgroup, then
$\Gamma $ belongs to $\mathcalc $.

\item \label{cal C and normal subgroups}
If $\Gamma$ contains a normal subgroup $\Delta$
which belongs to $\mathcalc$,
then $\Gamma $ belongs to $\mathcalc$.

\item
\label{amalgamated product}
Let $\Gamma$ be the amalgamated product
$\Gamma_0 \ast_{\Delta} \Gamma_1$ for a common subgroup
$\Delta$ of $\Gamma_0$ and $\Gamma_1$. Suppose that $\Delta$, $\Gamma_0$
and $\Gamma_1$ belong to ${\mathcal C}$ and that $\capac_0(\Delta) \le 0$,
then $\Gamma$ belongs to ${\mathcal C}$.

\item \label{cal C is L closed}
$L({\mathcal C} \cup {\mathcal F}) = {\mathcal C} \cup L{\mathcal F}$.

\item \label{groups with capac > 0}
Suppose for the group $\Gamma$ that $\capac_0(\Gamma) > 0$ or that
it contains $\zz^n$ for some $n \ge 1$ as a normal subgroup.
Then $\Gamma$ belongs to ${\mathcal C}$ and moreover
$$ \capac_p(\Gamma) \le 1 \hspace{15mm} \mbox{ also holds for } p \ge 1.$$

\end{enumerate}
\end{theorem}

\proof

\noindent
\ref{not in cal C}.) $H_0^{\Gamma}(E \Gamma ; \mathcaln (\Gamma) )$ is not cofinal-measurable
if $\Gamma$ is (locally) finite.

\noindent
\ref{cal C is L closed}.) We have to show for a group
$\Gamma \in L({\mathcal C} \cup {\mathcal F})$
which is not locally finite that $\Gamma$ belongs to ${\mathcal C}$.
Since any infinite finitely generated subgroup
$\Gamma'\subset \Gamma$ belongs to
${\mathcal C}$ by assumption and these groups
form a cofinal system of subgroups, we get from Theorem
\ref{estimates on capac, extensions and local information}.\ref{capac and finitely generated subgroups} that
$H_p^{\Gamma}(E \Gamma ; \mathcaln (\Gamma) )$ is cofinal-measurable.
and
\begin{equation*}
\capac_1(\Gamma) =  \liminf \{\capac_1(\Gamma')\} \leq 1 .
\end{equation*}
If $c_0(\Gamma )\leq 0$ we are done,
otherwise we know from \ref{zero-th capacity}.\ref{capac = 0^-} that $\Gamma$ is
locally virtually nilpotent but not locally finite and the result follows from
\ref{estimates on capac, extensions and local information}.\ref{capac_1 for virtually nilpotent groups}.

\noindent
\ref{cal C and normal subgroups}.) We get from
Theorem \ref{estimates on capac, extensions and local information}.\ref{capac
and normal subgroups}
that $H^{\Gamma}_p(E\Gamma;\mathcaln(\Gamma))$ is {\cofinalm} for $p \ge 0$ and
\begin{eqnarray}
\capac_1(\Gamma) & \le &\capac_0(\Delta) + \capac_1(\Delta) ~ \le ~ 1.
\label{eqn 3.7}
\end{eqnarray}
It remains to show
$\capac_0(\Gamma) + \capac_1(\Gamma) \le 1$.
If $c_0(\Gamma) \leq 0$, this follows from (\ref{eqn 3.7}).
If $\capac_0(\Gamma) > 0$ we are again in the case, where $\Gamma$ is locally virtually
nilpotent but not locally finite and can apply \ref{estimates on capac, extensions and local
  information}.\ref{capac_1 for virtually nilpotent groups}.

\noindent
\ref{cal C and cal E}.)
We use the following description of the
class of elementary-amenable groups
\cite[Lemma 3.1]{Kropholler-Linnell-Moody(1988)}. Let ${\mathcal B}$ be
the class of all groups which are finitely generated and virtually free abelian.
Define for each ordinal $\alpha$
\begin{eqnarray*}
{\mathcal E}_0 ~ =  ~ & \{1\}& ;
\\
{\mathcal E}_{\alpha} ~ = ~ & (L{\mathcal E}_{\alpha-1}){\mathcal B} &
\mbox{ , if $\alpha$ is a successor ordinal}
\\
{\mathcal E}_{\alpha} ~ = ~ &  \cup_{\beta < \alpha} {\mathcal E}_{\beta}&
\mbox{ , if $\alpha$ is a limit ordinal}.
\end{eqnarray*}
Then ${\mathcal E} = \cup_{\alpha \ge 0} {\mathcal E}_{\alpha}$
is the class of elementary
amenable groups.
For any elementary amenable group $\Gamma$ there is a least ordinal $\alpha$
with $\Gamma \in {\mathcal E}_{\alpha}$ and we use transfinite induction to show
that $\Gamma$ belongs to ${\mathcal C}$, provided that $\Gamma$ is infinite and
contains no
infinite locally finite subgroup.\par

The induction begin $\alpha = 0$ is obvious. If $\alpha$ is a limit ordinal,
the induction step is clear. Suppose that $\alpha$ is a successor ordinal.
Then there is an extension
$1 \longrightarrow \Delta \longrightarrow
\Gamma \longrightarrow \pi \longrightarrow 1$
such that $\Delta \in L{\mathcal E}_{\alpha -1}$ and $\pi \in {\mathcal B}$.
Every finitely generated subgroup $\Delta' \subset \Delta$ belongs to
${\mathcal E}_{\alpha -1}$. By assumption $\Gamma$ and therefore
$\Delta$ and $\Delta^{\prime}$ contain no infinite locally finite subgroup. The induction hypothesis implies that
$\Delta^{\prime}$ is finite or
belongs to ${\mathcal C}$. From assertions
\ref{cal C is L closed}.) we get
$\Delta \in L({\mathcal C} \cup {\mathcal F}) = {\mathcal C} \cup L{\mathcal F}$. Since $\Delta$ is not
infinite locally finite either $\Delta$ is finite or $\Delta \in
{\mathcal C}$. In the second case apply \ref{cal C and normal subgroups}.). If
 $\Delta$ is finite $\pi$ must be infinite. Since $\mathcal
 F\mathcal B=\mathcal B$ we have $\Gamma$ infinite virtually free
 abelian. In this case the statement follows from Theorem
 \ref{estimates on capac, extensions and local information}.\ref{capac and
   subgroups of finite index} and \ref{estimates on capac, extensions
   and local information}.\ref{capac_1 for virtually nilpotent
   groups}.

\noindent
\ref{amalgamated product}.)
This follows from the long Mayer-Vietoris sequence,
Theorem \ref{basic properties of capacity}.\ref{capacity
and exact sequences}, Lemma
\ref{elementary properties of special and
cofinal-special} and
Lemma \ref{capacity and induction}. One has to argue as above for the case $c_0 ( \Gamma ) >0$.

\noindent
\ref{groups with capac > 0}.)
Suppose $c_0( \Gamma) > 0$. Then $\Gamma$ is locally virtually nilpotent but not locally finite
and  the
statement is  \ref{estimates on capac, extensions and local information}.\ref{capac_1 for
  virtually nilpotent groups}.
Now suppose $\Gamma$ contains $\zz^n$ as a normal
subgroup for $n \geq 1$. Then the result follows from
\ref{estimates on
  capac, extensions and local information}.\ref{Z^n}  and \ref{estimates on
  capac, extensions and local information}.\ref{capac and normal
  subgroups}.
\qed


\typeout{--------------------- section 3 -----------------------}
\section{Final Remarks}
\label{Final Remarks}

\begin{remark} \label{possible improvement of Theorem}
We mention that the proof for Theorem
\ref{estimate on capac_1 of an
elementary amenable group}.\ref{cal C and cal E} goes also
through if one enlarges the class $\mathcal E$ as defined by transfinite
induction in the proof by substituting
the class ${\mathcal B}$ of virtually abelian groups
by any bigger class ${\mathcal B}'$ with the properties
that ${\mathcal B}' \subset {\mathcal C}$,
and ${\mathcal F}{\mathcal B}' = {\mathcal B}'$.
\qed
\end{remark}

\begin{remark} \label{exceptional groups}
Let $1 \longrightarrow \Delta \xrightarrow{i} \Gamma \xrightarrow{i} \zz
\longrightarrow 1$ be an extension of groups. Suppose
that $\Delta$ is locally finite. Then $H_p(E\Gamma;\mathcaln(\Gamma))$
is trivial for $p \ge 2$ and for $p=1$ equal to the kernel $K$
of the $\mathcaln(\Gamma)$-map
$$\mathcaln(\Gamma) \otimes_{\cc \Delta} \cc \longrightarrow
\mathcaln(\Gamma) \otimes_{\cc \Delta} \cc:\quad u \otimes n \mapsto
u(t-1) \otimes n$$
for some $t \in \mathcaln(\Gamma)$ which maps to a generator of $\zz$ under $p$.
If we would know that $K$ is trivial, then $\Gamma$ would belong to ${\mathcal C}$
and it would suffice in Theorem
\ref{estimate on capac_1 of an
elementary amenable group}.\ref{cal C and cal E}
to assume that $\Gamma$ itself
is not locally finite instead of assuming that
$\Gamma$ contains no infinite locally finite subgroup.
\qed
\end{remark}

\begin{remark} \label{capac_p(Gamma) le 1}
So far we know no counterexample to the following statement:

\noindent
If $\Gamma$ is elementary amenable  and not locally-finite,
then $H_p(E\Gamma;\mathcaln(\Gamma))$ is {\cofinalm} for all $p \ge 0$
and
\begin{eqnarray}
\capac_p(\Gamma) & \le & 1; \hspace{10mm} \mbox{ for } p \ge 0 .
\label{sjfsdgfasjdsgfas}
\end{eqnarray}

To prove this, it suffices to
show inequality \ref{sjfsdgfasjdsgfas} for any group $\Gamma$
such that there is an extension
$1 \longrightarrow \Delta \longrightarrow \Gamma
\longrightarrow \zz \longrightarrow 1$ with a group $\Delta$
which already satisfies inequality \ref{sjfsdgfasjdsgfas}.
Then the proof of Theorem
\ref{estimate on capac_1 of an elementary
amenable group}.\ref{cal C and cal E}
would go through.

On the other
hand one can construct a fiber bundle
$F \longrightarrow E \longrightarrow S^1$ of closed manifolds
with simply-connected fiber $F$ such that $\capac_p(\widetilde{E})$
is arbitrary large.
This follows from the observation that in this case one can read off
$\capac_p(\widetilde{E})$ from the automorphism of $H_p(F)$ induced by
the monodromy map $F \longrightarrow F$ and one can realize
any element in $GL(n,\zz)$ as the automorphism induced on
$H_3(\prod_{i=1}^n S^3)$ by an automorphism of $\prod_{i=1}^n S^3$.
\qed
\end{remark}




\typeout{--------------------- references -----------------------}


\end{document}